\documentclass[a4,11pt]{article}
\usepackage{amsmath}
\usepackage{graphics}
\usepackage{latexsym}
\usepackage{amsfonts}
\numberwithin{equation}{section}
\renewcommand{\baselinestretch}{1.40}
\title{
The Finite Horizon Optimal Multi-Modes Switching Problem: the
Viscosity Solution Approach.}
\author{Brahim EL-ASRI \thanks{Universit\'e du Maine, D\'ept. de Math\'ematiques, Equipe Stat. et Processus,
Avenue Olivier Messiaen, 72085 Le Mans, Cedex 9, France. e-mails:
brahim.el\_Asri.etu@univ-lemans.fr }\,\,\, \,and \, Said HAMADENE
\thanks{Universit\'e du Maine, D\'ept. de Math\'ematiques, Equipe
Stat. et Processus, Avenue Olivier Messiaen, 72085 Le Mans, Cedex 9,
France. e-mail: hamadene@univ-lemans.fr}}
\begin{document}
\date{}
\maketitle
\newtheorem{theo}{Theorem}
\newtheorem{problem}{Problem}
\newtheorem{pro}{Proposition}
\newtheorem{cor}{Corollary}
\newtheorem{axiom}{Definition}
\newtheorem{rem}{Remark}
\newtheorem{lem}{Lemma}
\newcommand{\brm}{\begin{rem}}
\newcommand{\erm}{\end{rem}}
\newcommand{\beth}{\begin{theo}}
\newcommand{\eeth}{\end{theo}}
\newcommand{\bl}{\begin{lem}}
\newcommand{\el}{\end{lem}}
\newcommand{\bp}{\begin{pro}}
\newcommand{\ep}{\end{pro}}
\newcommand{\bcor}{\begin{cor}}
\newcommand{\ecor}{\end{cor}}
\newcommand{\be}{\begin{equation}}
\newcommand{\ee}{\end{equation}}
\newcommand{\beq}{\begin{eqnarray*}}
\newcommand{\eeq}{\end{eqnarray*}}
\newcommand{\beqa}{\begin{eqnarray}}
\newcommand{\eeqa}{\end{eqnarray}}
\newcommand{\dg}{\displaystyle \delta}
\newcommand{\cm}{\cal M}
\newcommand{\cF}{{\cal F}}
\newcommand{\cR}{{\cal R}}
\newcommand{\bF}{{\bf F}}
\newcommand{\tg}{\displaystyle \theta}
\newcommand{\w}{\displaystyle \omega}
\newcommand{\W}{\displaystyle \Omega}
\newcommand{\vp}{\displaystyle \varphi}
\newcommand{\ig}[2]{\displaystyle \int_{#1}^{#2}}
\newcommand{\integ}[2]{\displaystyle \int_{#1}^{#2}}
\newcommand{\produit}[2]{\displaystyle \prod_{#1}^{#2}}
\newcommand{\somme}[2]{\displaystyle \sum_{#1}^{#2}}
\newlength{\inter}
\setlength{\inter}{\baselineskip} \setlength{\baselineskip}{7mm}
\newcommand{\no}{\noindent}
\newcommand{\rw}{\rightarrow}
\def \ind{1\!\!1}
\def \R{I\!\!R}
\def \cadlag {{c\`adl\`ag}~}
\def \esssup {\mbox{ess sup}}
\begin{abstract}

In this paper we show existence and uniqueness of a solution for a
system of $m$ variational partial differential inequalities with
inter-connected obstacles. This system is the deterministic version
of the Verification Theorem of the Markovian optimal $m$-states
switching problem. The switching cost functions are arbitrary. This
problem is in relation with the valuation of firms in a financial
market.
\end{abstract}

\no{\bf AMS Classification subjects}: 60G40 ; 62P20 ; 91B99 ; 91B28
; 35B37 ; 49L25.
\medskip

\no {$\bf Keywords$}: Real options; Backward stochastic differential
equations; Snell envelope; Stopping times ; Switching; Viscosity
solution of PDEs; Variational inequalities.

\section {Introduction}In this work we are concerned with the following system of
$m$ variational partial differential inequalities with
inter-connected obstacles: \be\label{sysintro} \left\{
\begin{array}{l}
\min\{v_i(t,x)- \max\limits_{j\in{\cal
I}^{-i}}(-g_{ij}(t,x)+v_j(t,x)),\\\qquad\qquad
-\partial_tv_i(t,x)-{\cal A}v_i(t,x)-
\psi_i(t,x)\}=0,\,\forall (t,x)\in [0,T]\times \R^k,\,\,i\in{\cal I}=\{1,...,m\},\\
v_i(T,x)=0.
\end{array}\right.
\ee where $g_{ij}$, $\psi_i$ are continuous functions, $\cal A$ an
infinitesimal generator associated with a diffusion process and
finally ${\cal I}^{-i}:=\{1,...,i-1,i+1,...,m\}$.

This system is the deterministic version of the Verification Theorem
of the optimal multi-modes switching problem in finite horizon. This
problem, of real option type, can be introduced with the help of the
following example:

Assume we have a power station/plant which produces electricity and
which has several modes of production, e.g., the lower, the middle
and the intensive modes. The price of electricity in the market,
given by an adapted stochastic process $(X_t)_{t \leq T}$,
fluctuates in reaction to many factors such as demand level, weather
conditions, unexpected outages and so on. On the other hand,
electricity is non-storable, once produced it should be almost
immediately consumed. Therefore, as a consequence, the station
produces electricity in its instantaneous most profitable mode known
that when the plant is in mode $i\in {\cal I}$, the yield per unit
time is given by means of $\psi_i$ and, on the other hand, switching
the plant from the mode $i$ to the mode $j$ is not free and
generates expenditures given by $g_{ij}$ and possibly by other
factors in the energy market. So the manager of the power plant
faces two main issues:

$(i)$ when should her decide to switch the production from its
current mode to another one?

$(ii)$ to which mode the production has to be switched when the
decision of switching is made?
\medskip

\noindent In other words she faces the issue of finding the optimal
strategy of management of the plant. This issue is in relation with
the price of the power plant in the energy market.
\medskip

For decades, optimal switching problems have attracted a lot of
research activity (see e.g. \cite{[BE], [BO1],[BS], [CL],[DX], [D],
[DP], [DH], [DHP], [DZ], [DZ2], [HJ], hib,ht, [LP],porchet,porchet2,
[RM], shi, [TY], tri1, tri, dz} and the references therein).
Especially in connection with valuation of firms and, questions
related to the structural profitability of investment project or an
industry whose production depends on the fluctuating market price of
a number of underlying commodities or assets, ... Several variants
of the problem we deal with here have been considered. In order to
tackle those problems, authors use mainly two approaches. Either a
probabilistic one (\cite{[DH], [DHP], [HJ], ht, porchet,porchet2})
or an   approach which uses partial differential inequalities (PDIs
for short) (\cite{[BE],[BO1],[CL],[DZ2],[LP], dz, [TY]}).

The PDIs approach turns out to study and to solve, in some sense,
the system of $m$ PDIs with inter-connected obstacles
(\ref{sysintro}). Amongst the papers which consider the same problem
as ours, and in the framework of viscosity solutions approach, the
most elaborated works are certainly the ones by Tang and Yong
\cite{[TY]}, on the one hand, and by Djehiche et al. \cite{[DHP]},
on the other hand. In \cite{[TY]}, the authors show existence and
uniqueness of a solution for (\ref{sysintro}). Nevertheless the
paper suffers from two facts: $(i)$ the growth exponent at infinity
of the functions $\psi_i$ should be smaller that $2$ ; $(ii)$ the
switching cost functions $g_{ij}$ should not depend on $x$. The
first issue of \cite{[TY]} has been treated by Djehiche et
al.\cite{[DHP]} since in their paper the authors show existence of
the solution for \ref{sysintro} in the case when the growth of the
functions $\psi_i$ is of arbitrary polynomial type. The second issue
of \cite{[TY]}, i.e. considering the case when $g_{ij}$ depending
also on $x$, was right now, according to our knowledge, an open
problem. Note that in \cite{[DHP]}, the question of uniqueness is
not addressed. Therefore the main objective of our work, and this is
the novelty of the paper, is to show existence and uniqueness of a
solution in viscosity sense for the system when the functions
$\psi_i$ and $g_{ij}$ are continuous depending also on $x$ and
satisfy an arbitrary polynomial growth condition. We show also that
the solution is unique in the class of continuous functions with
polynomial growth.

This paper is organized as follows:

In Section 2, we formulate the problem and we give the related
definitions.  In Section 3, we introduce the optimal switching
problem under consideration and give its probabilistic Verification
Theorem. It is expressed by means of a Snell envelope of processes.
Then we introduce the approximating scheme which enables to
construct a solution for the Verification Theorem. Moreover we give
some properties of that solution, especially the dynamic programming
principle. Section 4 is devoted to the connection between the
optimal switching problem, the Verification Theorem and the system
of PDIs (\ref{sysintro}). This connection is made through backward
stochastic differential equations with one reflecting obstacle in
the case when randomness comes from a solution of a standard
stochastic differential equation. Further we provide some estimate
for the optimal strategy of the switching problem which in
combination with the dynamic programming principle plays a crucial
role. Finally we show that system (\ref{sysintro}) has a solution.
In Section 5, we show that the solution of (\ref{sysintro}) is
unique in the class of continuous functions which satisfy a
polynomial growth condition. $\Box$

\section{Assumptions and formulation of the problem}
Throughout this paper $T$ (resp. $k$) is a fixed real (resp.
integer) positive constant. Let us now consider the followings:
\medskip

\indent $(i)$ $b:[0,T]\times\R^k\rightarrow \R^{k}$ and
$\sigma:[0,T]\times\R^k\rightarrow \R^{k\times d}$ are two
continuous functions for which there exists a constant $C\geq 0$
such that for any $t\in [0, T]$ and $x, x'\in \R^k$\be
\label{regbs1}|b(t,x)|+ |\sigma(t,x)|\leq C(1+|x|) \quad \mbox{ and
} \quad |\sigma(t,x)-\sigma(t,x')|+|b(t,x)-b(t,x')|\leq C|x-x'|\ee

$(ii)$ for $i,j \in {\cal I}=\{1,...,m\}$,
$g_{ij}:[0,T]\times\R^k\rightarrow \R$ and
$\psi_i:[0,T]\times\R^k\rightarrow \R$ are continuous functions and
of polynomial growth, $i.e.$, there exist some positive constants
$C$ and $\gamma$ such that for each $i,j\in \cal I$: \be
\label{polycond} |\psi_i(t,x)|+|g_{ij}(t,x)|\leq C(1+|x|^\gamma),\,
\, \forall (t,x)\in [0,T]\times \R^k. \ee Moreover we assume that
there exists a constant $\alpha >0$ such that for any $(t,x)\in
[0,T]\times \R^k$, \be  min\{g_{ij}(t,x), i,j \in {\cal I},\quad i
\neq j \} \geq \alpha .\ee
\medskip

We now consider the following system of $m$ variational inequalities
with inter-connected obstacles:  $\forall \,\,i\in {\cal I}$ \be
\label{sysvi1} \left\{
\begin{array}{l}
\min\left \{v_i(t,x)- \max\limits_{j\in{\cal I}^{-i}}(-g_{ij}(t,x)+v_j(t,x)),-\partial_tv_i(t,x)-
{\cal A}v_i(t,x)-\psi_i(t,x)\right\}=0,\\
v_i(T,x)=0,
\end{array}\right.
\ee where ${\cal I}^{-i}:={\cal I}-\{i\}$ and ${\cal A}$ is the
following infinitesimal generator:
\begin{equation}
\label{generateur} {\cal A}=\frac{1}{2}\sum_{i,j=1,k}(\sigma
\sigma^*)_{ij}(t,x)\frac{\partial^2} {\partial x_i \partial
x_j}+\sum_{i=1,k} b_i(t,x)\frac{\partial}{\partial x_i}\,;
\end{equation}
hereafter the superscript $(^*)$ stands for the transpose, $Tr$ is
the trace operator and finally $<x,y>$ is the inner product of
$x,y\in \R^k$.
\medskip

The main objective of this paper is to focus on the existence and
uniqueness of the solution in viscosity sense of (\ref{sysvi1})
whose definition is:

\begin{axiom} Let $(v_1,...,v_m)$ be a $m$-uplet of continuous functions defined on
$[0,T]\times \R^k$, $\R$-valued and such that such that $v_i(T,x)=0$
for any $x\in \R^k$ and $i\in {\cal I}$. The $m$-uplet
$(v_1,...,v_m)$ is called:
\begin{itemize}
\item [$(i)$] a viscosity supersolution (resp. subsolution) of the system (\ref{sysvi1})
if for each fixed $i\in {\cal I}$, for any $(t_0,x_0)\in [0,T]\times
\R^k$ and any function $\varphi_i \in C^{1,2}([0,T]\times \R^k)$
such that $\varphi_i(t_0,x_0)=v_i(t_0,x_0)$ and $(t_0,x_0)$ is a
local maximum of $\varphi_i -v_i$ (resp. minimum), we have:
\be
\begin{array}{l}
\min\left\{v_i(t_0,x_0)- \max\limits_{j\in{\cal
I}^{-i}}(-g_{ij}(t_0,x_0)+v_j(t_0,x_0)),\right.\\\qquad\qquad \qquad
\left.-\partial_t \varphi_i(t_0,x_0)-{\cal
A}\varphi_i(t_0,x_0)-\psi_i(t_0,x_0)\right\}\geq 0 \quad
(\mbox{resp.} \leq 0).
\end{array}
\ee
\item [$(ii)$] a viscosity solution if it is both a viscosity supersolution and
subsolution. $\Box$
\end{itemize}
\end{axiom}

There is an equivalent formulation of this definition (see e.g.
{\cite{[CIL]}) which we give because it will be useful later. So
firstly we define the notions of superjet and subjet of a continuous
function $v$.
\begin{axiom}
Let $v \in C((0,T)\times \R^k)$, $(t,x)$ an element of $(0,T)\times
\R^k$ and finally $S_k$ the set of $k \times k$ symmetric matrices.
We denote by $J^{2,+} v(t,x)$ (resp. $J^{2,-} v(t,x)$), the
superjets (resp. the subjets) of $v$ at $(t,x)$, the set of triples
$(p,q,X)\in \R\times \R^k \times S_k$ such that:
$$\begin{array}{c}
v(s,y)\leq v(t,x) + p(s-t)+\langle q,y-x \rangle +\frac{1}{2}\langle
X(y-x),y-x\rangle+o(|s-t|+|y-x|^2) \\
(resp.\quad v(s,y)\geq v(t,x) + p(s-t)+\langle q,y-x \rangle
+\frac{1}{2}\langle X(y-x),y-x\rangle+o(|s-t|+|y-x|^2)). \Box
\end{array}$$
\end{axiom}
Note that if $\varphi-v$ has a local maximum (resp. minimum) at
$(t,x)$, then we obviously have:
$$\left(D_t \varphi(t,x),D_x \varphi(t,x),D^{2}_{xx}\varphi(t,x)\right)
\in J^{2,-} v (t,x) \,\,\, (\mbox{resp. } J^{2,+} v (t,x)). \Box$$

We now give an equivalent definition of a viscosity solution of the
parabolic system with inter-connected obstacles (\ref{sysvi1}).
\begin{axiom}Let $(v_1,...,v_m)$ be a $m$-uplet of continuous functions defined
on $[0,T]\times \R^k$, $\R$-valued and such that
$(v_1,...,v_m)(T,x)=0$ for any $x\in \R^k$. The $m$-uplet
$(v_1,...,v_m)$ is called a viscosity supersolution (resp.
subsolution) of (\ref{sysvi1}) if for any $i\in {\cal I}$, $(t,x)\in
(0,T)\times \R^k$ and $(p,q,X)\in J^{2,-} v_i (t,x)$ (resp. $J^{2,+}
v_i (t,x)$),
$$ min \left\{v_i(t,x)-\max\limits_{j\in{\cal I}^{-i}}(-g_{ij}(t,x)+
v_j(t,x)),-p -\frac{1}{2}Tr[\sigma^{*} X \sigma] -\langle b,q
\rangle-\psi_{i}(t,x)\right\}\geq 0 \,\,(resp. \leq 0).$$ It is
called a viscosity solution it is both a viscosity subsolution and
supersolution .$\Box$
\end{axiom}

As pointed out previously we will show that system (\ref{sysvi1})
has a unique solution in viscosity sense. This system is the
deterministic version of the optimal $m$-states switching problem
which is well documented in \cite{[DHP]} and which we will describe
briefly in the next section.
\section{The  optimal $m$-states  switching problem}
\subsection{Setting of the problem}
Let $(\Omega, {\cal F}, P)$ be a fixed probability space on which is
defined a standard $d$-dimensional Brownian motion $B=(B_t)_{0\leq
t\leq T}$ whose natural filtration is $(\cF_t^0:=\sigma \{B_s, s\leq
t\})_{0\leq t\leq T}$. Let $ \bF=(\cF_t)_{0\leq t\leq T}$ be the
completed filtration of $(\cF_t^0)_{0\leq t\leq T}$ with the
$P$-null sets of ${\cal F}$, hence $(\cF_t)_{0\leq t\leq T}$
satisfies the usual conditions, $i.e.$, it is right continuous and
complete. Furthermore, let:

- ${\cal P}$ be the $\sigma$-algebra on $[0,T]\times \Omega$ of
$\bF$-progressively measurable sets;

- ${\cal M}^{2,k}$ be the set of $\cal P$-measurable and
$\R^k$-valued processes $w=(w_t)_{t\leq T}$ such that
$E[\int_0^T|w_s|^2ds]<\infty$  and ${\cal S}^2$  be the set of $\cal
P$-measurable, continuous processes ${w}=({w}_t)_{t\leq T}$ such
that $E[\sup_{t\leq T}|{w}_t|^2]<\infty$;

-  for any stopping time $\tau \in [0,T]$, ${\cal T}_\tau$ denotes
the set of all stopping times $\theta$ such that $\tau \leq \theta
\leq T$. $\Box$
\medskip

The problem of multiple switching can be described through an
example as follows. Assume  we have a plant which produces a
commodity, $e.g.$ a power station which produces electricity. Let
${\cal I} $ be the set of all possible activity modes of the
production of the commodity. A management strategy of the plant
consists, on the one hand, of the choice of a sequence of
nondecreasing stopping times $(\tau_n)_{n\geq1}$ $(i.e. \tau_n \leq
\tau_{n+1}$ and $\tau_0 = 0)$ where the manager decides to switch
the activity from its current mode to another one. On the other
hand, it consists of the choice of the mode $\xi_n$, a r.v. ${\cal
F}_{\tau_n}$-measurable with values in ${\cal I}$, to which the
production is switched at $\tau_n$ from its current mode. Therefore
the admissible management strategies of the plant are the pairs
$(\delta,\xi):=((\tau_n)_{n\geq 1},(\xi_n)_{n\geq 1})$ and the set
of these strategies is denoted by $\cal D$.

Let now $X:=(X_t)_{0\leq t\leq T}$ be an $\cal P$-measurable,
$\R^k$-valued continuous stochastic process which stands for the
market price of $k$ factors which determine the market price of the
commodity. On the other hand, assuming that the production activity
is in mode 1 at the initial time $t = 0$, let $(u_t)_{t\leq T}$
denote the indicator of the production activity's mode at time $t\in
[0, T]$ :
\begin{equation}
u_t=\ind_{[0,\tau_1]}(t)+\sum_{n\geq1}\xi_n
\ind_{(\tau_{n},\tau_{n+1}]}(t).
\end{equation}
Then for any $t\leq T$, the state of the whole economic system
related to the project at time $t$ is represented by the vector :
\begin{equation}
\begin{array}{ll}
(t, X_t, u_t)\in [0,T]\times \R^k \times {\cal I}.
\end{array}
\end{equation}

Finally, let $\psi_i(t,X_t)$ be the instantaneous profit when the
system is in state $(t, X_t, i)$, and for $i,j \in {\cal I} \quad
i\neq j$, let $g_{ij}(t,X_t)$ denote the switching cost of the
production at time $t$ from current mode $i$ to another mode $j$.
Then if the plant is run under the strategy
$(\delta,\xi)=((\tau_n)_{n\geq 1},(\xi_n)_{n\geq 1})$ the expected
total profit is given by:
$$\begin{array}{l} J(\delta,\xi)=E[\integ{0}{T}\psi_{u_s}(s,X_s)ds
-\sum_{n\geq 1} g_{u_{\tau_{n-1}}u_{\tau_n}}(\tau_{n},X_{\tau_{n}}) \ind_{[\tau_{n}<T]}].
\end{array}
$$
Therefore the problem we are interested in is to find an optimal
strategy, $i.e$, a strategy $(\delta^*,\xi^*)$ such that
$J(\delta^*,\xi^*)\ge J(\delta,\xi)$ for any $(\delta,\xi)\in \cal
D$.
\medskip

Note that in order that the quantity $J(\delta,\xi)$ makes sense we
assume throughout this paper that for any $i,j\in {\cal I}$ the
processes $(\psi_i(t,X_t))_{t\leq T}$ and $(g_{ij}(t,X_t))_{t\leq
T}$ belong to ${\cal M}^{2,1}$ and ${\cal S}^{2}$ respectively. On
the other hand there is a bijective correspondence between the pairs
$(\delta,\xi)$ and the pairs $(\delta,u)$. Therefore throughout this
paper one refers indifferently to $(\delta,\xi)$ or $(\delta,u)$.
\subsection{ The Verification Theorem}

To tackle the problem described above Djehiche et al. \cite{[DHP]}
have introduced a Verification Theorem which is expressed by means
of Snell envelope of processes. The Snell envelope of a stochastic
process $(\eta_t)_{t\leq T}$ of ${\cal S}^2$ (with a possible
positive jump at $T$) is the lowest supermartingale
$R(\eta):=(R(\eta)_t)_{t\leq T}$ of ${\cal S}^2$ such that for any
$t\leq T$, $R(\eta)_t\geq \eta_t$. It has the following expression:
$$\forall t\leq T, R(\eta)_t=esssup_{\tau \in {\cal
T}_t}E[\eta_\tau|\bF_t] \mbox{ and satisfies }R(\eta)_T=\eta_T.$$
For more details on the Snell envelope notion on can see e.g.
\cite{[CK], Elka, ham}.
\medskip

The Verification Theorem for the $m$-states optimal switching
problem is the following:
\begin{theo}(\cite{[DHP]}, Th.1)\label{thmverif}
\noindent Assume that there exist $m$ processes $(Y^i:=(Y^i_t)_{0\le
t\leq T}, i=1,...,m)$  of ${\cal S}^2$ such that: \be
\begin{array}{l}
\label{eqvt} \forall t\leq T,\,\,Y^i_t=\esssup_{\tau \geq
t}E[\int_t^\tau\psi_i(s,X_s)ds +\max\limits_{j\in {\cal
I}^{-i}}(-g_{ij}(\tau,X_\tau)+Y^j_\tau)1_{[\tau <T]}|\cF_t],\quad
Y^i_T=0.
\end{array}
\ee Then:
\begin{itemize}
\item[$(i)$]
$ Y^1_0=\sup \limits_{(\delta,\xi) \in {\cal D}}J(\delta,u). $
\item[$(ii)$] Define the sequence of $\bF$-stopping times $\delta^{*}=(\tau_n^*)_{n\geq 1}$ as follows :
$$
\begin{array}{lll}
\tau^*_{1}&=&\inf\{s\geq 0,\quad Y^1_s=\max\limits_{j\in{{\cal I}^{-1}}}(-g_{1j}(s,X_s)+Y^j_s)\}\wedge T,\\
\tau^*_{n}&=&\inf\{s\geq \tau^*_{n-1},\quad Y^{u_{\tau^*_{n-1}}}_s=\max\limits_{k\in {\cal I}\backslash \{u_{\tau^*_{n-1}}\}}(-g_{u_{\tau^*_{n-1}}k}(s,X_s)+Y^k_s)\}\wedge T,
\quad \mbox{for}\quad n\geq 2,
\end{array}
$$
where:
\begin{itemize}
\item[$\bullet$] $u_{\tau^*_1}=\sum \limits_{j\in {\cal I}} j \ind_ {\{\max\limits_{k\in {\cal I}^{-1}} (-g_{1k}(\tau^*_1,X_{\tau^*_1})+Y^k_{\tau^*_1})=-g_{1j}(\tau^*_1,X_{\tau^*_1})+Y^j_{\tau^*_1}\}};$
\item[$\bullet$] for any $n\geq1$ and $t\geq \tau^*_n,$
$Y^{u_{\tau^*_n}}_t =\sum\limits_{j\in {\cal I}} \ind_{[u_{\tau^*_n}=j]}Y^j _t$
\item[$\bullet$] for any $n\geq 2, \,\,u_{\tau^*_n}=l$ on the set
$$\left\{\max\limits_{k\in {\cal I}
\backslash \{{ u_{\tau^*_{n-1}}}\}} (-g_{u_{\tau^*_{n-1}}
k}(\tau^*_{n},X_{\tau^*_{n}})+ Y^k_{\tau^*_{n}})=-g_{u_{\tau^*_{n-1}
l}}(\tau^*_n,X_{\tau^*_n})+Y^l_{\tau^*_n}\right\}$$with \,
$g_{u_{\tau^*_{n-1} k}}(\tau^*_n,X_{\tau^*_n})=
\sum\limits_{j\in {\cal I}}\ind_{[u_{\tau^*_{n-1}}=j]}g_{j k}(\tau^*_n,X_{\tau^*_n})$ and ${\cal I}\backslash \{u_{\tau^*_{n-1}}\}=\sum\limits_{j\in {\cal I}}\ind_{[u_{\tau^*_{n-1}}=j]}{\cal I}^{-j}$.
\end{itemize}
Then the strategy $(\delta^*,u^*)$ is optimal i.e.
$J(\delta^*,u^*)\geq J(\delta,u)$ for any $(\delta,u)\in \cal D$.
$\Box$
\end{itemize}
\end{theo}
\medskip

The issue of existence of the processes $Y^1,...,Y^m$ which satisfy
(\ref{eqvt}) is also addressed in \cite{[DHP]}. Also for $n\geq 0$
let us define the processes $(Y^{1,n},...,Y^{m,n})$ recursively as
follows: for $i\in {\cal I}$ we set,
\begin{equation}\label{y0}
Y^{i,0}_t=\esssup_{\tau\geq t}E[\integ{t}{\tau}\psi_i(s,X_s)ds|{\cal
F}_t],\,\,0\le t\leq T,
\end{equation}
and for $n\geq 1$,
\begin{equation}
\label{eq24} Y^{i,n}_t=\mbox{ess sup}_{\tau\geq t}
E[\integ{t}{\tau}\psi_i(s,X_s)ds+\max\limits_{k\in {\cal
I}^{-i}}(-g_{ik}(\tau,X_\tau)+Y^{k,n-1}_\tau) \ind_{[\tau<T]}|{\cal
F}_t],\,\,0\le t\leq T.
\end{equation}
Then the sequence of processes $((Y^{1,n},...,Y^{m,n}))_{n\geq 0}$
have the following properties:
\begin{pro} (\cite{[DHP]}, Pro.3 and Th.2)
\begin{itemize}
\item[$(i)$] for any $i\in {\cal I}$ and
$n\geq 0$, the processes $Y^{1,n},...,Y^{m,n}$ are well-posed,
continuous and belong to ${\cal S}^2$, and verify
\be\label{croi}\forall t\le T,\,\,Y^{i,n}_t\leq Y^{i,n+1}_t\leq
E[\int_t^T\{\max_{i=1,m}|\psi_i(s,X_s)|\}ds|{\cal F}_t];\ee

\item[$(ii)$] there exist $m$ processes $Y^1,...,Y^m$ of ${\cal S}^2$ such
that for any $i\in {\cal I}$:
\begin{itemize}
\item[$(a)$]  $\forall t\leq T$, $Y^i_t=\lim_{n\rightarrow
\infty}\nearrow Y^{i,n}_t$ and
$$
E[\sup_{s\leq T}|Y^{i,n}_s-Y^i_s|^2]\to 0 \quad\mbox{ as }\quad n\to
+\infty $$
\item[$(b)$] $\forall t\leq T$,
\begin{equation}
\label{eq26}{Y}^{i}_t=\mbox{ess sup}_{\tau\geq
t}E[\integ{t}{\tau}\psi_i(s,X_s)ds+ \max\limits_{k\in {\cal
I}^{-i}}(-g_{ik}(\tau,X_\tau)+{Y}^{k}_\tau) \ind_{[\tau<T]}|{\cal
F}_t]
\end{equation}
i.e. ${Y}^1,...,{Y}^m$ satisfy the Verification Theorem
\ref{thmverif} ;
\item[$(c)$] $\forall t\leq T$,
\begin{equation}
\label{eq27}{Y}^{i}_t=esssup_{(\delta,\xi)\in {\cal
D}^i_t}E[\integ{t}{T}\psi_{u_s}(s,X_s)ds -\sum_{n\geq 1}
g_{u_{\tau_{n-1}}u_{\tau_n}}(\tau_{n},X_{\tau_{n}})
\ind_{[\tau_{n}<T]}]|{\cal F}_t]
\end{equation}
where ${\cal D}^i_t=\{(\delta,\xi)=((\tau_n)_{n\geq
1},(\xi_n)_{n\geq 1}) \mbox { such that } u_0 =i \mbox{ and }
\tau_1\geq t \}$. This characterization means that if at time $t$
the production activity is in its regime $i$ then the optimal
expected profit is $Y^i_t$.
\item[$(d)$] the processes $Y^1,...,Y^m$ verify the dynamical programming principle of the
$m$-states optimal switching problem, $i.e.$, $\forall t\leq T$,
\begin{equation}\!\!\!\!\!
\label{dpp}
\begin{array}{l}
Y^i_t=\esssup_{(\delta,u) \in {\cal D}_t^i}
E[\integ{t}{\tau_n}\psi_{u_s}(s,X_s)ds -\sum_{1\leq k \leq n}
g_{u_{\tau_{k -1}}u_{{\tau_k}}}({\tau}_{k},X_{{\tau}_{k}})
\ind_{[\tau_{k}<T]}+\ind_{[\tau_n <T]}Y^{u_{\tau_n}}_{\tau_n}|{\cal
F}_t].\Box
\end{array}
\end{equation}
\end{itemize}
\end{itemize}
\end{pro}

Note that except $(ii-d)$, the proofs of the other points are given
in \cite{[DHP]}. The proof of $(ii-d)$ can be easily deduced in
using relation (\ref{eq26}). Actually from (\ref{eq26}) for any
$i\in {\cal I}$, $t\in [0,T]$ and $(\delta,\xi)\in {\cal D}^i_t$ we
have: \be \label{eq28}Y^i_t\geq
E[\integ{t}{\tau_n}\psi_{u_s}(s,X_s)ds -\sum_{1\leq k \leq n}
g_{u_{\tau_{k -1}}u_{u_{\tau_k}}}({\tau}_{k},X_{{\tau}_{k}})
\ind_{[\tau_{k}<T]}+\ind_{[\tau_n <T]}Y^{u_{\tau_n}}_{\tau_n}|{\cal
F}_t]. \ee Next using the optimal strategy we obtain the equality
instead of inequality in (\ref{eq28}). Therefore the relation
(\ref{dpp}) holds true. $\Box$
\begin{rem} \label{unic}Note that the characterization (\ref{eq27}) implies that the processes $Y^1,...,Y^m$ of ${\cal S}^2$ which satisfy the Verification
Theorem are unique.
\end{rem}
\section{Existence of a solution for the system of variational inequalities}
\subsection{Connection with BSDEs with  one reflecting barrier}
Let $(t,x)\in [0,T]\times \R^k$ and let $(X^{tx}_s)_{s\leq T}$ be
the solution of the following standard SDE:
\begin{equation}\label{sde}
dX^{tx}_s=b(s,X_s^{tx})ds+\sigma(s,X_s^{tx})dB_s \mbox{ for }t\leq
s\leq T\mbox{ and }X_s^{tx}=x \mbox{ for }s\leq t\end{equation}where
the functions $b$ and $\sigma$ are the ones of (\ref{regbs1}). These
properties of $\sigma$ and $b$ imply in particular that the process
$(X^{tx}_s)_{0\le s\leq T}$ solution of the standard SDE (\ref{sde})
exists and is unique, for any $t\in [0, T]$ and $x\in \R^k$.

The operator $\cal A$  that is appearing in (\ref{generateur}) is
the infinitesimal generator associated with $X^{t,x}$. In the
following result we collect some properties of $X^{t,x}$.

\bp \label{estimx} (see e.g. \cite{[RY]}) The process $X^{tx}$ satisfies
the following estimates:
\begin{itemize}
\item [$(i)$] For any $q\geq 2$, there exists a constant $C$ such that
\begin{equation}\label{estimat1}
E[\sup_{0\le s\leq T}|X^{tx}_s|^q]\leq C(1+|x|^q).
\end{equation}
\item[$(ii)$] There exists a constant $C$ such that for any $t,t'\in [0,T]$ and $x,x'\in \R^k$,
\begin{equation}\label{estimat2}
E[\sup_{0\le s\leq T}|X^{tx}_s-X^{t'x'}_s|^2]\leq
C(1+|x|^2)(|x-x'|^2+|t-t'|). \Box
\end{equation}
\end{itemize}
\ep

We are going now to introduce the notion of a BSDE with one
reflecting barrier introduced in \cite{[EKal]}. This notion will
allow us to make the connection between the variational inequalities
system (\ref{sysvi1}) and the $m$-states optimal switching problem
described in the previous section.
\medskip

So let us introduce the deterministic functions $f:[0,T]\times
\R^{k+1+d}\rightarrow \R$, $h:[0,T]\times \R^{k}\rightarrow \R$ and
$g:\R^{k}\rightarrow \R$ continuous, of polynomial growth and such
that $h(x,T)\leq g(x)$. Moreover we assume that for any $(t,x)\in
[0,T]\times \R^{k}$, the mapping $(y,z)\in \R^{1+d}\mapsto
f(t,x,y,z)$ is uniformly Lipschitz. Then we have the following
result related to BSDEs with one reflecting barrier:
\begin{theo}(\cite{[EKal]}, Th. 5.2 and 8.5) For any $(t,x)\in [0,T]\times \R^k$, there exits a unique triple
of processes $(Y^{t,x},Z^{t,x},K^{t,x})$ such that:
\begin{equation}\left\{
\begin{array}{l}
Y^{tx}, K^{tx}\in {\cal S}^2 \mbox{ and }Z^{tx}\in {\cal
M}^{2,d};\,K^{tx}
\mbox{ is  non-decreasing and }K^{tx}_0=0,\\
Y^{tx}_s=g(X^{tx}_T)+\int_{s}^{T}f(r,X_r^{tx},Y^{tx}_r,Z^{tx}_r)dr-\int_{s}^{T}Z^{tx}_rdB_r+
K_T^{tx}-K^{tx}_s, \,\, s\le T\\
Y^{tx}_s\geq h(s,X^{tx}_s),\, \forall s\le T\mbox{ and }
\int_{0}^{T}(Y^{tx}_r-h(r,X^{tx}_r))dK^{tx}_r=0.
\end{array}
\right. \end{equation} Moreover the following characterization of
$Y^{t,x}$ as a Snell envelope holds true: \be
\label{snellenv}\forall s\leq T,\,\,Y^{t,x}_s=esssup_{\tau \in {\cal
T}_t}E[\int_t^\tau f(r,X_r^{tx},Y^{tx}_r,Z^{tx}_r)dr
+h(\tau,X^{tx}_\tau)\ind_{[\tau<T]}+g(X^{tx}_T)\ind_{[\tau=T]}|{\cal
F}_s].\ee

On the other hand there exists a deterministic continuous with
polynomial growth function $u: [0,T]\times \R^{k}\rightarrow \R$
such that:$$\forall s\in [t,T], Y^{t,x}_s=u(s,X^{t,x}_s).$$ Moreover
the function $u$ is the unique viscosity solution in the class of
continuous function with polynomial growth of the following PDE with
obstacle:$$ \left\{
\begin{array}{l}
\min\{u(t,x)- h(t,x), -\partial_tu(t,x)-{\cal
A}u(t,x)-f(t,x,u(t,x),\sigma (t,x)^*\nabla u(t,x))\}=0,
\\u(T,x)=g(x).\Box
\end{array}\right.
$$
\end{theo}
\subsection{Existence of a solution for the system of variational inequalities}

Let $(Y^{1,tx}_s,...,Y^{m,tx}_s)_{0\le s\leq T}$ be the processes
which satisfy the Verification Theorem \ref{thmverif} in the case
when the process $X\equiv X^{t,x}$. Therefore using the
characterization (\ref{snellenv}), there exist processes $K^{i,tx}$
and $Z^{i,tx}$, $i\in {\cal I}$, such that the triples ($Y^{i,tx},
Z^{i,tx},K^{i,tx})$ are unique solutions (thanks to Remark
\ref{unic}) of the following reflected BSDEs: for any $i=1,...,m$ we
have,
\begin{equation}\left\{
\begin{array}{l}
Y^{i,tx}, K^{i,tx}\in {\cal S}^2 \mbox{ and }Z^{i,tx}\in {\cal
M}^{2,d};\,K^{i,tx}
\mbox{ is  non-decreasing and }K^{i,tx}_0=0,\\
Y^{i,tx}_s=\int_{s}^{T}\psi_i(r,X_r^{tx})du-\int_{s}^{T}Z^{i,tx}_rdB_r+K_T^{i,tx}-K^{i,tx}_s, \,\,\, 0\le s\le T,\,\,Y^{i,tx}_T=0,\\
Y^{i,tx}_s\geq \max\limits_{j\in{\cal I}^{-i}}(-g_{ij}(s,X_s^{tx})+Y^{j,tx}_s),\,\,0\le s\le T,\\
\int_{0}^{T}(Y^{i,tx}_r-\max\limits_{j\in{\cal
I}^{-i}}(-g_{ij}(r,X_r^{tx})+Y^{j,tx}_r))dK^{i,tx}_r=0.
\end{array}
\right. \end{equation} Moreover we have the following result.
\begin{pro}There are deterministic functions $v^1,...,v^m$
$:[0,T]\times \R^k\rightarrow \R$ such that:
$$\forall (t,x)\in [0,T]\times \R^k, \forall s\in [t,T],
Y_s^{i,tx}=v^i(s,X^{tx}_s), \,\,i=1,...,m.$$ Moreover the functions
$v^i$, $i=1,...,m,$ are lower semi-continuous and of polynomial
growth.
\end{pro}
$Proof$: For $n\geq 0$ let $(Y^{n,1,tx}_s,...,Y^{n,m,tx}_s)_{0\le
s\leq T}$ be the processes constructed in (\ref{y0})-(\ref{eq24}).
Therefore using an induction argument and Theorem 2 there exist
deterministic continuous with polynomial growth functions $v^{i,n}$
($i=1,...,m$) such that for any $(t,x)\in [0,T]\times \R^k$,
$\forall s\in [t,T]$, $Y^{n,i,tx}_s=v^{i,n}(s,X^{tx}_s)$. Using now
inequality (\ref{croi}) we get:
$$Y^{n,i,tx}_t\le Y^{n+1,i,tx}_t\leq
E[\int_t^T\{\max_{i=1,m}|\psi_i(s,X^{tx}_s)|\}ds]$$ since
$Y^{n,i,tx}_t$ is deterministic. Therefore combining the polynomial
growth of $\psi_i$ and estimate (\ref{estimat1}) for $X^{tx}$ we
obtain:
$$v^{i,n}(t,x)\leq v^{i,n+1}(t,x)\leq C(1+|x|^p)$$for some constants
$C$ and $p$ independent of $n$. In order to complete the proof it is
enough now to set $v^i(t,x):=\lim_{n\rightarrow \infty}v^{i,n}(t,x),
(t,x)\in [0,T]\times \R^k$ since $Y^{i,n,tx}\nearrow Y^{i,tx}$ as $n
\rightarrow \infty$. $\Box$
\medskip

We are now going to focus on the continuity of the functions
$v^1,...,v^m$. But first let us deal with some properties of the
optimal strategy which exist thanks to Theorem 1.

\bp \label{optimal-s} Let $(\delta,u)=((\tau_n)_{n\geq
1},(\xi_n)_{n\geq 1})$ be an optimal strategy, then there exist two
positive constant $C$ and $p$ which does not depend on $t$ and $x$
such that: \be \label{estiopt}\forall n\geq 1,\,\,P[\tau_{n}<T]\leq
\frac {C(1+|x|^p)}{n}.\ee \ep $Proof$: Recall the characterization
of (\ref{eq27}) that reads as:
$$
\begin{array}{l}
Y^{1,tx}_0=sup_{(\delta,u) \in {\cal
D}}E[\int_0^T\psi_{u_r}(r,X_r^{tx})dr-\sum_{k\geq 1}
g_{u_{\tau_{k-1}}u_{\tau_k}}(\tau_{k},X^{tx}_{\tau_{k}})
\ind_{[\tau_{k}<T]}].
\end{array}
$$
Now if $(\delta,u)=((\tau_n)_{n\geq 1},(\xi_n)_{n\geq 1})$ is the
optimal strategy then we have: $$
\begin{array}{l}
Y^{1,tx}_0=E[\int_0^T\psi_{u_r}(r,X_r^{tx})dr-\sum_{k\geq 1}
g_{u_{\tau_{k-1}}u_{\tau_k}}(\tau_{k},X^{tx}_{\tau_{k}})
\ind_{[\tau_{k}<T]}].
\end{array}
$$
Taking into account that $g_{ij}\geq \alpha >0$ for any $i\neq j$ we
obtain:
$$\begin{array}{ll}
E[\sum_{k=1,n} \alpha \ind_{[\tau_{k}<T]}]+ Y^{1,tx}_0&\leq
E[\int_0^T\psi_{u_r}(r,X_r^{tx})dr-\sum_{k\geq n+1}
g_{u_{\tau_{k-1}}u_{\tau_k}}(\tau_{k},X^{tx}_{\tau_{k}})
\ind_{[\tau_{k}<T]}].
\end{array}
$$
But for any $k\le n$, $[\tau_{n}<T]\subset[\tau_{k}<T]$
then:$$\begin{array}{ll} \alpha n P[\tau_{n}<T]+ Y^{1,tx}_0&\leq
E[\int_0^T\psi_{u_r}(r,X_r^{tx})dr-\sum_{k\geq n+1}
g_{u_{\tau_{k-1}}u_{\tau_k}}(\tau_{k},X^{tx}_{\tau_{k}})
\ind_{[\tau_{k}<T]}]\\{}&\leq E[\int_0^T\psi_{u_r}(r,X_r^{tx})dr].
\end{array}
$$and then
$$
\begin{array}{ll}
n\alpha P[\tau_{n}<T] &\leq E[\int_0^T\mid
\psi_{u_r}(r,X_r^{tx})\mid dr]-Y^{1,tx}_0\\{}&\leq  E[\int_0^T\mid
\psi_{u_r}(r,X_r^{tx})\mid dr]-Y^{1,0,tx}_0.
\end {array}
$$
Finally taking into account the facts that $\psi_i$ and $Y^{1,0,tx}$
are of polynomial growth and estimate (\ref{estimat1}) for $X^{tx}$ to
obtain the desired result. Note that the polynomial growth of
$Y^{1,0,tx}$ stems from Proposition 3. $\Box$ \brm The estimate
(\ref{estiopt}) is also valid for the optimal strategy if at the
initial time the state of plant is an arbitrary $i\in {\cal I }$.
$\Box$ \erm
\medskip

Next for $i\in {\cal I}$ let
$(y^{i,tx}_s,z^{i,tx}_s,k^{i,tx}_s)_{s\leq T}$ be the processes
defined as follows:
\begin{equation}\label{y2} \left\{
\begin{array}{l}
y^{i,tx}, k^{i,tx}\in {\cal S}^2 \mbox{ and }z^{i,tx}\in {\cal
M}^{2,d};\,k^{i,tx}
\mbox{ is  non-decreasing and }k^{i,tx}_0=0,\\
y^{i,tx}_s=\integ{s}{T}\psi_i(r,X_r^{tx})\ind_{[r\geq t]}dr-\integ{s}{T}z^{i,tx}_rdB_r+k_T^{i,tx}-k^{i,tx}_s, \,\,\,
0\le s\le T, \,\, y^{i,tx}_T=0\\ y^{i,tx}_s\geq l^{i,tx}_{s}:=
\max\limits_{j\in{\cal I}^{-i}}\{-g_{ij}(t\vee s,X_{t\vee
s}^{tx})+y^{j,tx}_s)\},\,\,\,
\forall s\le T, \\
\integ{0}{T}(y^{i,tx}_r-l^{i,tx}_{r})dk^{i,tx}_r=0.
\end{array}
\right. \end{equation} The existence of
$(y^{i,tx},z^{i,tx},k^{i,tx}), i\in {\cal I}$, is obtained in the
same way as the one of $(Y^{i,tx},Z^{i,tx},K^{i,tx})$. On the other
hand, thanks to uniqueness (see once more Remark \ref{unic}), for
any $(t,x)\in [0,T]\times \R^k$, for any $s\in [0,t]$ we have
$y^{i,tx}_s=Y^{i,tx}_t$, $z^{i,tx}_s=0$ and $k^{i,tx}_s=0$.
\medskip

We are now ready to give the main Theorem of this article. \beth The
functions $(v^1,...,v^m):[0,T]\times \R^k\rightarrow \R$ are
continuous and solution in viscosity sense of the system of
variational inequalities with inter-connected obstacles
(\ref{sysvi1}).\eeth $Proof$: First let us focus on continuity and
let us show that $v^1$ is continuous. The same proof will be valid
for the continuity of the other functions $v^i$ ($i=2,...,m$).

First the characterization (\ref{eq27}) implies that:
$$
y^{1,tx}_0=\sup_{(\delta,\xi)\in {\cal
D}}E[\int_0^T\psi_{u_s}(s,X^{tx}_s)\ind_{[s\geq t]}ds -\sum_{n\geq
1} g_{u_{\tau_{n-1}}u_{\tau_n}}(\tau_{n}\vee t,X^{tx}_{\tau_{n}\vee
t}) \ind_{[\tau_{n}<T]}]
$$
On the other hand an optimal strategy $(\delta^*,\xi^*)$ exists and
satisfies the estimates (\ref{estiopt}) with the same constants $C$
and $p$. Next let $\epsilon >0$ and $(t',x')\in B((t,x),\epsilon)$
and let us consider the following set of strategies:$$ \tilde
D:=\{(\delta,\xi)=((\tau_n)_{n\geq 1}, (\xi_n)_{n\geq 0}) \in {\cal
D} \mbox{ such that } \forall n\geq 1, P[\tau_n<T] \leq
\frac{C(1+(\epsilon +|x|)^p)}{n}\}.$$ Therefore the strategy
$(\delta^*,\xi^*)$ belongs to $\tilde D$ and then we have:

$$\begin{array}{ll} y^{1,tx}_0&=\sup_{(\delta,\xi)\in {\tilde
D}}E[\int_0^T\psi_{u_s}(s,X^{tx}_s)\ind_{[s\geq t]}ds -\sum_{n\geq
1} g_{u_{\tau_{n-1}}u_{\tau_n}}(\tau_{n}\vee t,X^{tx}_{\tau_{n}\vee
t}) \ind_{[\tau_{n}<T]}]\\
{}&=sup_{(\delta,u) \in {\tilde D}}
E[\int_{0}^{\tau_n}\psi_{u_s}(s,X^{tx}_s)\ind_{[s\geq t]}ds
\\{}&\qquad\qquad\qquad-\sum_{1\leq k \leq n} g_{u_{\tau_{k
-1}}u_{\tau_k}}(t\vee{\tau}_{k} ,X^{tx}_{{t\vee \tau}_{k}})
\ind_{[\tau_{k}<T]}+\ind_{[\tau_n
<T]}y^{u_{\tau_n},tx}_{\tau_n}]\end{array}
$$
and
$$\begin{array}{ll}
y^{1,t'x'}_0&=\sup_{(\delta,\xi)\in {\tilde
D}}E[\int_0^T\psi_{u_s}(s,X^{t'x'}_s)\ind_{[s\geq t']}ds
-\sum_{n\geq 1} g_{u_{\tau_{n-1}}u_{\tau_n}}(\tau_{n}\vee
t',X^{t'x'}_{\tau_{n}\vee t'}) \ind_{[\tau_{n}<T]}]\\
{}&=sup_{(\delta,u) \in {\tilde D}}
E[\int_{0}^{\tau_n}\psi_{u_s}(s,X^{t'x'}_s)\ind_{[s\geq t']}ds
\\{}&\qquad\qquad\qquad-\sum_{1\leq k \leq n} g_{u_{\tau_{k
-1}}u_{\tau_k}}(t'\vee{\tau}_{k} ,X^{t'x'}_{{t'\vee \tau}_{k}})
\ind_{[\tau_{k}<T]}+\ind_{[\tau_n
<T]}y^{u_{\tau_n},t'x'}_{\tau_n}].\end{array}
$$
The second equalities follow from the dynamical programming
principle. It follows that:\be \label{eqcont}
\begin{array}{ll}y^{1,t'x'}_0-y^{1,tx}_0&\leq
sup_{(\delta,u) \in {\tilde D}}
E[\int_{0}^{\tau_n}\{\psi_{u_s}(s,X^{t'x'}_s)\ind_{[s\geq
t']}-\psi_{u_s}(s,X^{tx}_s)\ind_{[s\geq t]}\}ds\\
{}&\qquad -\sum_{1\leq k \leq n} \{g_{u_{\tau_{k
-1}}u_{\tau_k}}(t'\vee{\tau}_{k} ,X^{t'x'}_{{t'\vee \tau}_{k}})
-g_{u_{\tau_{k -1}}u_{\tau_k}}(t\vee{\tau}_{k} ,X^{tx}_{{t\vee
\tau}_{k}})\} \ind_{[\tau_{k}<T]}\\{}&\qquad +\ind_{[\tau_n
<T]}\{y^{u_{\tau_n},t'x'}_{\tau_n}-y^{u_{\tau_n},tx}_{\tau_n}\}]\end{array}\ee
Next w.l.o.g we assume that $t'<t$. Then from (\ref{eqcont}) we
deduce that: \be \label{eqcont2}
\begin{array}{ll}y^{1,t'x'}_0-y^{1,tx}_0&\leq
sup_{(\delta,u) \in {\tilde D}}
E[\int_{0}^{\tau_n}\{\psi_{u_s}(s,X^{t'x'}_s)-\psi_{u_s}(s,X^{tx}_s))\ind_{[s\geq
t]}+\psi_{u_s}(s,X^{t'x'}_s)\ind_{[t'\leq s< t]}\}ds\\
{}&\qquad -\sum_{1\leq k \leq n} \{g_{u_{\tau_{k
-1}}u_{\tau_k}}(t'\vee{\tau}_{k} ,X^{t'x'}_{{t'\vee \tau}_{k}})
-g_{u_{\tau_{k -1}}u_{\tau_k}}(t\vee{\tau}_{k} ,X^{tx}_{{t\vee
\tau}_{k}})\} \ind_{[\tau_{k}<T]}\\{}&\qquad +\ind_{[\tau_n
<T]}\{y^{u_{\tau_n},t'x'}_{\tau_n}-y^{u_{\tau_n},tx}_{\tau_n}\}]\\
{}&\leq
E[\int_{0}^{T}\max_{j=1,m}|\psi_{j}(s,X^{t'x'}_s)-\psi_{j}(s,X^{tx}_s))|\ind_{[s\geq
t]}+\max_{j=1,m}|\psi_j(s,X^{t'x'}_s)|\ind_{[t'\leq s< t]}\}ds\\
{}&\qquad +n\max_{i\neq j\in {\cal I}}\{\sup_{s\leq T}|g_{ij}(t'\vee
s ,X^{t'x'}_{t'\vee s }) -g_{ij}(t\vee s ,X^{tx}_{t\vee s})|\}
\\{}&\qquad +sup_{(\delta,u) \in {\tilde D}}(P[\tau_n
<T])^{\frac{1}{2}}(2E[(y^{u_{\tau_n},t'x'}_{\tau_n})^2+(y^{u_{\tau_n},tx}_{\tau_n})^2])^{\frac{1}{2}}.
\end{array}
\ee In the right-hand side of (\ref{eqcont2}) the first term
converges to $0$ as $(t',x')\rightarrow (t,x)$. Next let us show
that for any $i\neq j \in {\cal I}$, $$E[\sup_{s\leq
T}|g_{ij}(t'\vee s ,X^{t'x'}_{t'\vee s }) -g_{ij}(t\vee s
,X^{tx}_{t\vee s})|]\rightarrow 0 \mbox{ as }(t',x')\rw (t,x).$$
Actually for any $\varrho>0$ we have: $$\begin{array}{l}
|g_{ij}(t'\vee s ,X^{t'x'}_{t'\vee s }) -g_{ij}(t\vee s
,X^{tx}_{t\vee s})|\leq
|g_{ij}(t'\vee s ,X^{t'x'}_{t'\vee s }) -g_{ij}(t\vee s
,X^{t'x'}_{t'\vee s})|\ind_{[|X^{t'x'}_{t'\vee s }|\leq \varrho]}+\\
\qquad\qquad  |g_{ij}(t'\vee s ,X^{t'x'}_{t'\vee s }) -g_{ij}(t\vee
s ,X^{t'x'}_{t'\vee s})|\ind_{[|X^{t'x'}_{t'\vee s }|\geq
\varrho]}+|g_{ij}(t\vee s ,X^{t'x'}_{t'\vee s }) -g_{ij}(t\vee s
,X^{tx}_{t\vee s})|.
\end{array}
$$
Therefore we have:
$$\begin{array}{l}
E[\sup_{s\leq T}|g_{ij}(t'\vee s ,X^{t'x'}_{t'\vee s })
-g_{ij}(t\vee s ,X^{tx}_{t\vee s})|]\leq\\ \qquad E[\sup_{s\leq T}\{
|g_{ij}(t'\vee s ,X^{t'x'}_{t'\vee s }) -g_{ij}(t\vee s
,X^{t'x'}_{t'\vee s})|\ind_{[|X^{t'x'}_{t'\vee s }|\leq \varrho]}\}]+\\
\qquad  E[\sup_{s\leq T}\{|g_{ij}(t'\vee s ,X^{t'x'}_{t'\vee s })
-g_{ij}(t\vee s ,X^{t'x'}_{t'\vee
s})|\}\ind_{[\sup_{s\leq T}|X^{t'x'}_{t'\vee s }|\geq \varrho]}]+\\
\qquad E[\sup_{s\leq T}\{|g_{ij}(t\vee s ,X^{t'x'}_{t'\vee s })
-g_{ij}(t\vee s ,X^{tx}_{t\vee s})|\}\ind_{[\sup_{s\leq
T}|X^{t'x'}_{t'\vee s }|+\sup_{s\leq T}|X^{tx}_{t\vee s }|\geq
\varrho]}]+\\\qquad E[\sup_{s\leq T}\{|g_{ij}(t\vee s
,X^{t'x'}_{t'\vee s }) -g_{ij}(t\vee s ,X^{tx}_{t\vee
s})|\}\ind_{[\sup_{s\leq T}|X^{t'x'}_{t'\vee s }|+\sup_{s\leq
T}|X^{tx}_{t\vee s }|\leq \varrho]}]
\end{array}
$$
But since $g_{ij}$ is continuous then it is uniformly continuous on
$[0,T]\times \{x\in \R^k, |x|\leq \varrho\}.$ Henceforth for any
$\epsilon_1>0$ there exists $\eta_{\epsilon_1}>0$ such that for any
$|t-t'|<\eta_{\epsilon_1}$ we have:\be\label{terme1}\sup_{s\leq T}\{
|g_{ij}(t'\vee s ,X^{t'x'}_{t'\vee s }) -g_{ij}(t\vee s
,X^{t'x'}_{t'\vee s})|\ind_{[|X^{t'x'}_{t'\vee s }|\leq
\varrho]}\}\leq \epsilon_1.\ee Next using Cauchy-Schwarz's
inequality and then Markov's one with the second term we obtain: \be
E[\sup_{s\leq T}\{|g_{ij}(t'\vee s ,X^{t'x'}_{t'\vee s })
-g_{ij}(t\vee s ,X^{t'x'}_{t'\vee s})|\}\ind_{[\sup_{s\leq
T}|X^{t'x'}_{t'\vee s }|\geq \varrho]}]\leq
C(1+|x'|^p)\varrho^{-\frac{1}{2}} \ee where $C$ and $p$ are real
constants which are bound to the polynomial growth of $g_{ij}$ and
estimate (\ref{estimat1}). In the same way we have: \be
E[\sup_{s\leq T}\{|g_{ij}(t\vee s ,X^{t'x'}_{t'\vee s })
-g_{ij}(t\vee s ,X^{tx}_{t\vee s})|\}\ind_{[\sup_{s\leq
T}|X^{t'x'}_{t'\vee s }|+\sup_{s\leq T}|X^{tx}_{t\vee s }|\geq
\varrho]}]\leq C(1+|x|^p+|x'|^p) \varrho^{-\frac{1}{2}}\ee Finally
using the uniform continuity of $g_{ij}$ on compact subsets, the
continuity property (\ref{estimat2}) and the Lebesgue dominated
convergence theorem to obtain that \be \label{term4}E[\sup_{s\leq
T}\{|g_{ij}(t\vee s ,X^{t'x'}_{t'\vee s }) -g_{ij}(t\vee s
,X^{tx}_{t\vee s})|\ind_{[\sup_{s\leq T}|X^{t'x'}_{t'\vee s
}|+\sup_{s\leq T}|X^{tx}_{t\vee s }|\leq \varrho]}]\rw 0 \mbox{ as
}(t',x')\rw (t,x).\ee Taking now into account
(\ref{terme1})-(\ref{term4}) we have:
$$\limsup_{(t',x')\rw (t,x)}E[\sup_{s\leq T}|g_{ij}(t'\vee s ,X^{t'x'}_{t'\vee s })
-g_{ij}(t\vee s ,X^{tx}_{t\vee s})|]\leq
\epsilon_1+C(1+|x|^p)\varrho^{-\frac{1}{2}}.$$ As $\epsilon_1$ and
$\varrho$ are arbitrary then making $\epsilon_1\rightarrow 0$ and
$\varrho \rw +\infty$ to obtain that:
$$\lim_{(t',x')\rw (t,x)}E[\sup_{s\leq T}|g_{ij}(t'\vee s ,X^{t'x'}_{t'\vee s })
-g_{ij}(t\vee s ,X^{tx}_{t\vee s})|]=0.$$ Thus the claim is proved.

Finally let us focus on the last term in (\ref{eqcont2}). Since
$(\delta,u)\in \tilde D$ then:
$$\begin{array}{ll}sup_{(\delta,u) \in {\tilde D}}(P[\tau_n
<T])^{\frac{1}{2}}(2E[(y^{u_{\tau_n},t'x'}_{\tau_n})^2+(y^{u_{\tau_n},tx}_{\tau_n})^2])^{\frac{1}{2}}&\leq
n^{-\frac{1}{2}}\sup_{(\delta,u) \in {\tilde
D}}(2E[(y^{u_{\tau_n},t'x'}_{\tau_n})^2+(y^{u_{\tau_n},tx}_{\tau_n})^2])^{\frac{1}{2}}\\
{}&\leq Cn^{-\frac{1}{2}}(1+|x|^p+|x'|^p)\end{array}$$ where $C$ and
$p$ are appropriate constants which come from the polynomial growth
of $\psi_i$, $i\in {\cal I}$, estimate (\ref{estimat1}) for the
process $X^{tx}$ and inequality (\ref{croi}). Going back now to
(\ref{eqcont2}), taking the limit as $(t',x')\rw (t,x)$ to obtain:
$$
\limsup_{(t',x')\rw (t,x)}y^{1,t'x'}_0\leq y^{1,tx}_0
+Cn^{-\frac{1}{2}}(1+2|x|^p).$$As $n$ is arbitrary then putting
$n\rw +\infty$ to obtain:
$$
\limsup_{(t',x')\rw (t,x)}y^{1,t'x'}_0\leq y^{1,tx}_0.$$ It implies
that: $$\limsup_{(t',x')\rw
(t,x)}y^{1,t'x'}_0=Y^{1,t'x'}_{t'}=v^1(t',x')\leq
y^{1,tx}_0=Y^{1,tx}_{t}=v^1(t,x).
$$
Therefore $v^1$ is upper semi-continuous. But $v^1$ is also lower
semi-continuous, therefore it is continuous. In the same way we can
show that $v^2$,...,$v^m$ are continuous. As they are of polynomial
growth then taking into account Theorem 2 to obtain that
$(v^1,\dots,v^m)$ is a viscosity solution for the system of
variational inequalities with inter-connected obstacles
(\ref{sysvi1}). $\Box$
\section{Uniqueness of the solution of the system} We are going now to address the
question of uniqueness of the viscosity solution of the system
(\ref{sysvi1}). We have the following:

\beth \label{uni}The solution in viscosity sense of the system of
variational inequalities with inter-connected obstacles
(\ref{sysvi1}) is unique in the space of continuous functions on
$[0,T]\times R^k$ which satisfy a polynomial growth condition, i.e.,
in the space
$$\begin{array}{l}{\cal C}:=\{\varphi: [0,T]\times \R^k\rightarrow
\R, \mbox{ continuous and for any }\\\qquad \qquad\qquad(t,x), \,
|\varphi(t,x)|\leq C(1+|x|^\gamma) \mbox{ for some constants } C
\mbox{ and }\gamma\}.\end{array}$$ \eeth {\it Proof}. We will show
by contradiction that if $u_1,...,u_m$ and $w_1,...,w_m$ are a
subsolution and a supersolution respectively for (\ref{sysvi1}) then
for any $i=1,...,m$, $u_i\leq w_i$. Therefore if we have two
solutions of (\ref{sysvi1}) then they are obviously equal. Actually
for some $R>0$ suppose there exists
$(\overline{t},\overline{x},\overline{i})\in(0,T)\times B_R\times
{\cal I}$ $(B_R := \{x\in \R^k; |x|<R\})$ such that:
\begin{equation}
\label{comp-equ}
u_{\overline{i}}(\overline{t},\overline{x})-w_{\overline{i}}(\overline{t},\overline{x})=\eta>0.
\end{equation}Let us take
$\theta,\lambda$ and $\beta \in (0,1]$ small enough, so that the
following holds:
\be
\left\{
\begin{array}{llll}
\beta T<\frac{\eta}{5}\\
2\theta|\overline{x}|^{2\gamma + 2}< \frac{\eta}{5}\\
-\lambda w_{\overline{i}}(\overline{t},\overline{x})< \frac{\eta}{5}\\
\frac{\lambda}{\overline{t}}<\frac{\eta}{5}.
\end{array}
\right. \ee Here $\gamma$ is the growth exponent of the functions
which w.l.o.g we assume integer and $\geq 2$. Then, for a small
$\epsilon>0$, let us define:
\begin{equation}
\label{phi}
\Phi^i_{\epsilon}(t,x,y)=u_{i}(t,x)-(1-\lambda)w_{i}(t,y)-\frac{1}{2\epsilon}|x-y|^{2\gamma}
-\theta(|x|^{2\gamma + 2}+|y|^{2\gamma + 2})+\beta t -
\frac{\lambda}{t}.
\end{equation}
By the growth assumption on $u_i$ and $w_i$, there exists a $(t_{0},x_{0},y_{0},i_0)\in (0,T]\times \overline{B}_R \times \overline{B}_R \times {\cal I}$, such that:
$$\Phi^{i_0}_{\epsilon}(t_{0},x_{0},y_{0})=\max\limits_{(t,x,y,i)}\Phi^i_{\epsilon}(t,x,y).$$
On the other hand, from $2\Phi^{i_0}_{\epsilon}(t_{0},x_{0},y_{0})\geq \Phi^{i_0}_{\epsilon}(t_0,x_0,x_0)+\Phi^{i_0}_{\epsilon}(t_0,y_0,y_0)$, we have
\begin{equation}
\frac{1}{2\epsilon}|x_0 -y_0|^{2\gamma} \leq (u_{i_0}(t_0,x_0)-u_{i_0}(t_0,y_0))+(1-\lambda)(w_{i_0}(t_0,x_0)-w_{i_0}(t_0,y_0)),
\end{equation}
and consequently $\frac{1}{2\epsilon}|x_0 -y_0|^{2\gamma}$ is bounded, and as $\epsilon\rightarrow 0$, $|x_0 -y_0|\rightarrow 0$. Since $u_{i_0}$
and $w_{i_0}$ are uniformly continuous on $[0,T]\times \overline{B}_R$, then $\frac{1}{2\epsilon}|x_0 -y_0|^{2\gamma}\rightarrow 0$ as $\epsilon\rightarrow 0.$\\

Next let us show that $t_0 <T.$ Actually if $t_0 =T$ then,
$$
\Phi^{\overline{i}}_{\epsilon}(\overline{t},\overline{x},\overline{x})\leq \Phi^{i_0}_{\epsilon}(T,x_{0},y_{0}),$$
and,
$$
u_{\overline{i}}(\overline{t},\overline{x})-(1-\lambda)w_{\overline{i}}(\overline{t},\overline{x})-2\theta|\overline{x}|^{2\gamma
+ 2}+\beta \overline{t}- \frac{\lambda}{\overline{t}}\leq  \beta T -
\frac{\lambda}{T},
$$
since $u_{i_0}(T,x_0)=w_{i_0}(T,y_0)=0.$ Then thanks to
(\ref{comp-equ}) we have,
$$
\begin{array}{ll}
\eta &\leq -\lambda w_{\overline{i}}(\overline{t},\overline{x})+\beta T +2\theta|\overline{x}|^{2\gamma + 2}+\frac{\lambda}{\overline{t}}\\
\eta &< \frac{4}{5}\eta .
\end{array}
$$
which yields a contradiction and we have $t_0 \in (0,T)$. We now
claim that:
\begin{equation}
\label{visco-comp1} u_{i_0}(t_0,x_0)- \max\limits_{j\in{\cal
I}^{-i_0}}\{-g_{i_0 j}(t_0,x_0)+u_j(t_0,x_0)\} > 0.
\end{equation}
Indeed if
$$u_{i_0}(t_0,x_0)- \max\limits_{j\in{\cal I}^{-i_0}}\{-g_{i_0 j}(t_0,x_0)+u_j(t_0,x_0)\} \leq
0$$ then there exists $k \in {\cal I}^{-i_0}$ such that:
$$u_{i_0}(t_0,x_0) \leq -g_{i_0
k}(t_0,x_0)+u_k(t_0,x_0).$$
From the supersolution property of $w_{i_0}(t_0,y_0)$, we have
$$ w_{i_0}(t_0,y_0)\geq \max\limits_{j\in{\cal I}^{-i_0}}(-g_{i_0 j}(t_0,y_0)+w_j(t_0,y_0)) $$
then
$$ w_{i_0}(t_0,y_0)\geq -g_{i_0 k}(t_0,y_0)+w_k(t_0,y_0).$$
It follows that:
$$u_{i_0}(t_0,x_0)- (1-\lambda)w_{i_0}(t_0,y_0) -(u_{k}(t_0,x_0)-(1-\lambda)w_{k}(t_0,y_0))\leq (1-\lambda)g_{i_0 k}(t_0,y_0)-g_{i_0 k}(t_0,x_0).$$
Now since $g_{ij}\geq \alpha >0$, for every $i\neq j$, and taking into account of (\ref{phi}) to obtain:
$$\Phi^{i_0}_{\epsilon}(t_{0},x_{0},y_{0})-\Phi^{k}_{\epsilon}(t_{0},x_{0},y_{0})< -\alpha \lambda
+g_{i_0 k}(t_0,y_0)-g_{i_0 k}(t_0,x_0).$$
But this contradicts the definition of $i_0$, since $g_{i_0 k}$ is
uniformly continuous on $[0,T]\times \overline{B}_R$ and the claim
(\ref{visco-comp1}) holds.

Next let us denote
\begin{equation}
\varphi_{\epsilon}(t,x,y)=\frac{1}{2\epsilon}|x-y|^{2\gamma}+\theta(|x|^{2\gamma
+ 2}+|y|^{2\gamma + 2}) -\beta t + \frac{\lambda}{t}.
\end{equation}
Then we have: \be
\left\{
\begin{array}{lllll}\label{derive}
D_{t}\varphi_{\epsilon}(t,x,y)=-\beta- \frac{\lambda}{t^2},\\
D_{x}\varphi_{\epsilon}(t,x,y)= \frac{\gamma}{\epsilon}(x-y)|x-y|^{2\gamma-2} +\theta(2\gamma + 2)
x|x|^{2\gamma}, \\
D_{y}\varphi_{\epsilon}(t,x,y)= -\frac{\gamma}{\epsilon}(x-y)|x-y|^{2\gamma-2} +
\theta(2\gamma + 2)y|y|^{2\gamma},\\\\
B(t,x,y)=D_{x,y}^{2}\varphi_{\epsilon}(t,x,y)=\frac{1}{\epsilon}
\begin{pmatrix}
a_1(x,y)&-a_1(x,y) \\
-a_1(x,y)&a_1(x,y)
\end{pmatrix}+ \begin{pmatrix}
a_2(x)&0 \\
0&a_2(y)
\end{pmatrix} \\\\
\mbox{ with } a_1(x,y)=\gamma|x-y|^{2\gamma-2}I+\gamma(2\gamma -2)(x-y)(x-y)^* |x-y|^{2\gamma-4} \mbox{ and }\\
a_2(x)=\theta(2\gamma + 2)|x|^{2\gamma}I+2\theta \gamma(2\gamma +
2)xx^* |x|^{2\gamma-2 }.
\end{array}
\right. \ee Taking into account (\ref{visco-comp1}) then applying
the result by Crandall et al. (Theorem 8.3, {\cite{[CIL]}) to the
function $$
u_{i_0}(t,x)-(1-\lambda)w_{i_0}(t,y)-\varphi_{\epsilon}(t,x,y) $$ at
the point $(t_0,x_0,y_0)$, for any $\epsilon_1 >0$, we can find
$c,d\in \R$ and $X,Y \in S_k$, such that:

\be \label{lemmeishii}
\left\{
\begin{array}{lllll}
(c,\frac{\gamma}{\epsilon}(x_0-y_0)|x_0-y_0|^{2\gamma-2} +\theta(2\gamma + 2)x_0|x_0|^{2\gamma},X)
\in J^{2,+}(u_{i_0}(t_0,x_0)),\\
(-d,\frac{\gamma}{\epsilon}(x_0-y_0)|x_0-y_0|^{2\gamma-2} -\theta(2\gamma + 2)y_0|y_0|^{2\gamma },Y)\in J^{2,-}
((1-\lambda)w_{i_0}(t_0,y_0)),\\
c+d=D_{t}\varphi_{\epsilon}(t_0,x_0,y_0)=-\beta- \frac{\lambda}{t_0^2} \mbox{ and finally }\\
-(\frac{1}{\epsilon_1}+||B(t_0,x_0,y_0)||)I\leq
\begin{pmatrix}
X&0 \\
0&-Y
\end{pmatrix}\leq B(t_0,x_0,y_0)+\epsilon_1 B(t_0,x_0,y_0)^2.
\end{array}
\right. \ee Taking now into account (\ref{visco-comp1}), and the
definition of viscosity solution, we get:
$$\begin{array}{l}-c-\frac{1}{2}Tr[\sigma^{*}(t_0,x_0)X\sigma(t_0,x_0)]-\langle\frac{\gamma}{\epsilon}(x_0-y_0)|x_0-y_0|^{2\gamma-2}
+\\\qquad\qquad\qquad\qquad\qquad\theta(2\gamma +
2)x_0|x_0|^{2\gamma },b(t_0,x_0)\rangle-\psi_{i_0}(t_0,x_0)\leq 0
\mbox{ and
}\\d-\frac{1}{2}Tr[\sigma^{*}(t_0,y_0)Y\sigma(t_0,y_0)]-\langle\frac{\gamma}{\epsilon}(x_0-y_0)|x_0-y_0|^{2\gamma-2}
-\\\qquad\qquad\qquad\qquad\qquad\theta(2\gamma +
2)y_0|y_0|^{2\gamma},b(t_0,y_0)\rangle-(1-\lambda)\psi_{i_0}(t_0,y_0)\geq
0\end{array}$$ which implies that:
\begin{equation}
\begin{array}{llll}
\label{viscder}
-c-d&\leq \frac{1}{2}Tr[\sigma^{*}(t_0,x_0)X\sigma(t_0,x_0)-\sigma^{*}(t_0,y_0)Y\sigma(t_0,y_0)]\\
&\qquad +
\langle\frac{\gamma}{\epsilon}(x_0-y_0)|x_0-y_0|^{2\gamma-2},b(t_0,x_0)-b(t_0,y_0)\rangle\\&\qquad+\langle
\theta(2\gamma + 2)x_0|x_0|^{2\gamma },b(t_0,x_0)\rangle +\langle
\theta(2\gamma + 2)y_0|y_0|^{2\gamma },b(t_0,y_0)\rangle
\\&\qquad+\psi_i(t_0,x_0)-(1-\lambda)\psi_i(t_0,y_0).
\end{array}
\end{equation}
But from (\ref{derive}) there exist two constants $C$ and $C_1$ such
that:
$$||a_1(x_0,y_0)||\leq C|x_0 - y_0|^{2\gamma -2} \mbox{ and }(||a_2(x_0)||\vee ||a_2(y_0)||)\leq C_1 \theta.$$
As
$$B= B(t_0,x_0,y_0)= \frac{1}{\epsilon}
\begin{pmatrix}
a_1(x_{0},y_{0})&-a_1(x_{0},y_{0}) \\
-a_1(x_{0},y_{0})&a_1(x_{0},y_{0})
\end{pmatrix}+ \begin{pmatrix}
a_2(x_0)&0 \\
0&a_2(y_0)
\end{pmatrix}$$
then
$$B\leq \frac{C}{\epsilon}|x_0 - y_0|^{2\gamma -2}
\begin{pmatrix}
I&-I \\
-I&I
\end{pmatrix}+ C_1 \theta I.$$
It follows that:
\begin{equation}
B+\epsilon_1 B^2 \leq C(\frac{1}{\epsilon}|x_0 - y_0|^{2\gamma -2}+
\frac{\epsilon_1}{\epsilon^2}|x_0 - y_0|^{4\gamma -4})\begin{pmatrix}
I&-I \\
-I&I
\end{pmatrix}+ C_1\theta I
\end{equation}
where $C$ and $C_1$ which hereafter may change from line to line.
Choosing now $\epsilon_1=\epsilon$, yields the relation
\begin{equation}
\label{ineg_matreciel}
B+\epsilon_1 B^2 \leq \frac{C}{\epsilon}(|x_0 - y_0|^{2\gamma -2}+|x_0 - y_0|^{4\gamma -4})\begin{pmatrix}
I&-I \\
-I&I
\end{pmatrix}+ C_1\theta I.
\end{equation}
Now, from (\ref{regbs1}), (\ref{lemmeishii}) and
(\ref{ineg_matreciel}) we get:
$$\frac{1}{2}Tr[\sigma^{*}(t_0,x_0)X\sigma(t_0,x_0)-\sigma^{*}(t_0,y_0)
Y\sigma(t_0,y_0)]\leq \frac{C}{\epsilon}(|x_0 - y_0|^{2\gamma}+|x_0
- y_0|^{4\gamma -2}) +C_1 \theta(1+|x_0|^2+|y_0|^2).$$ Next $$
\langle\frac{\gamma}{\epsilon}(x_0-y_0)|x_0-y_0|^{2\gamma-2},b(t_0,x_0)-b(t_0,y_0)\rangle
\leq \frac{C^2}{\epsilon}|x_0 - y_0|^{2\gamma}$$ and finally,
$$\langle \theta(2\gamma + 2)x_0|x_0|^{2\gamma},b(t_0,x_0)\rangle +
\langle \theta(2\gamma + 2)y_0|y_0|^{2\gamma },b(t_0,y_0)\rangle
\leq \theta C(1+|x_0|^{2\gamma + 2}+|y_0|^{2\gamma + 2}).$$ So that
by plugging into (\ref{viscder}) and note that $\lambda >0$ we
obtain:
$$\begin{array}{l}\beta \leq \frac{C}{\epsilon}(|x_0 - y_0|^{2\gamma}+|x_0 - y_0|^{4\gamma -2}) +
C_1 \theta (1+|x_0|^2+|y_0|^2)+\frac{C^2}{\epsilon}|x_0 -
y_0|^{2\gamma}+ \\\qquad\qquad \theta C(1+|x_0|^{2\gamma +
2}+|y_0|^{2\gamma + 2})+
\psi_{i_0}(t_0,x_0)-(1-\lambda)\psi_{i_0}(t_0,y_0).\end{array}$$ By
sending $\epsilon\rightarrow0$, $\lambda\rightarrow0$, $\theta
\rightarrow0$ and taking into account of the continuity of
$\psi_{i_0}$ and $\gamma \geq 2$, we obtain $\beta \leq 0$ which is
a contradiction. The proof of Theorem \ref{uni} is now complete.
$\Box$
\medskip

As a by-product we have the following corollary: \bcor Let
$(v^1,...,v^m)$ be a viscosity solution of (\ref{sysvi1}) which
satisfies a polynomial growth condition then for $i=1,...,m$ and
$(t,x)\in [0,T]\times \R^k$, $$ v^i(t,x)= \sup_{(\delta,\xi)\in
{\cal D}^i_t}E[\integ{t}{T}\psi_{u_s}(s,X^{tx}_s)ds -\sum_{n\geq 1}
g_{u_{\tau_{n-1}}u_{\tau_n}}(\tau_{n},X^{tx}_{\tau_{n}})
\ind_{[\tau_{n}<T]}].
$$
\ecor

\no {\bf Acknowledgement}: the authors thank gratefully Prof.
J.Zhang for the fructuous discussions during the preparation of this
paper.$\Box$

\end{document}